# A real-time reactive framework for the surgical case sequencing problem

Belinda Spratt[*1] and Erhan Kozan[1]

**Abstract:** In this paper, we address the multiple operating room (**OR**) surgical case sequencing problem (**SCSP**). The objective is to maximise total OR utilisation during standard opening hours. This work uses a case study of a large Australian public hospital with long surgical waiting lists and high levels of non-elective demand. Due to the complexity of the SCSP and the size of the instances considered herein, heuristic techniques are required to solve the problem. We present constructive heuristics based on both a modified block scheduling policy and an open scheduling policy. A number of real-time reactive strategies are presented that can be used to maintain schedule feasibility in the case of disruptions. Results of computational experiments show that this approach maintains schedule feasibility in real-time, whilst increasing operating theatre (**OT**) utilisation and throughput, and reducing the waiting time of non-elective patients. The framework presented here is applicable to the real-life scheduling of OT departments, and we provide recommendations regarding implementation of the approach.

**Acknowledgements:** This research was funded by the Australian Research Council (**ARC**) Linkage Grant LP 140100394. Computational resources and services used in this work were provided by the High Performance Computing and Research Support Group, Queensland University of Technology, Brisbane, Australia.

---

[1] School of Mathematical Sciences, Queensland University of Technology (QUT)
2 George St, Brisbane, QLD 4000, Australia
* corresponding author, b.spratt@qut.edu.au  +617 3138 3035
Belinda Spratt ORCID: orcid.org/0000-0002-2522-2967
Erhan Kozan ORCID: orcid.org/0000-0002-3208-702X

# 1. Introduction

The surgical department is a highly uncertain environment. Whilst robust surgical schedules can reduce the likelihood of disruption, administrative staff require strategies to reschedule in the case that schedule deviations occur. These deviations include changes to staff availability (sick days, late arrivals), surgical durations, non-elective arrivals, and patient, staff and operating room (**OR**) based cancellations.

The approach presented here sequences patients, and assigns expected start and end times to their surgery. We allocate a qualified surgeon to each patient's surgery, while respecting surgeon availability. Whilst a single objective is considered explicitly, multiple performance measures are observed, in-line with hospital priorities (c.f. Cardoen et al. (2010)). The main objective is to maximise total time utilised during standard operating theatre (**OT**)[2] opening hours; however, we also consider total surgeon overtime, the average non-elective patient waiting time, and the number of additional elective patients scheduled. This makes the model a useful scheduling tool for administrative staff members as they can investigate schedules produced prioritising staff, resources, or patients, and various combinations of each.

Throughout the literature, when solving the real-time reactive surgical case sequencing problem (**SCSP**) surgeons and ORs are pooled resources. Our approach is novel as it incorporates availability and suitability of both ORs and surgeons, treating them as individual entities. We contribute further to existing literature by incorporating short-notice scheduling of elective patients where appropriate. In doing so, we are able to add previously unplanned elective surgeries to a schedule to reflect the short notice scheduling of transplant patients. The model presented in Section 3 is a bespoke formulation equivalent to a resource-constrained parallel-machine scheduling problem with identical machines, machine eligibility restrictions, and machine release dates. These are innovative additions to the current reactive scheduling literature.

---

[2] Note: the OT is the set of ORs.

We solve the model using constructive heuristics, based on an open scheduling policy and a block scheduling policy, which are well suited to the real-life surgical department. Thus, we make recommendations to hospital administrative staff members regarding ad-hoc changes to surgical schedules.

The layout of the paper is as follows. We review relevant OT planning and scheduling literature in Section 2, including the incorporation of stochasticity to create robust and reactive scheduling techniques. We provide a mixed integer programming (**MIP**) formulation of the real-time reactive SCSP in Section 3, based on the aforementioned case study. As the SCSP is NP-hard (Cardoen et al. , 2009) we present two constructive heuristics, and a set of real-time reactive rescheduling strategies in Section 4. We summarise the surgical department at the case study hospital in Section 5, including an overview of data and parameters (cf. Section 5.1). We examine computational results in Section 6, provide implementation recommendations in Section 7, and present concluding remarks in Section 8.

## 2. Literature Review

The administrative practices of surgical departments can have a large impact on hospital costs, patient outcomes, and the overall efficiency of a hospital. As such, a large number of papers exist which investigate and analyse some of the most significant issues in OT planning and scheduling. For a thorough review of the existing literature see Ferrand et al. (2014), Van Riet and Demeulemeester (2015), or Samudra et al. (2016). OT planning and scheduling problems can exist at the strategic, tactical, or operational levels. In this review, we focus on reactive and dynamic approaches to Surgical Case Sequencing (**SCS**) at the operational level of OT planning and scheduling.

One of the most significant challenges in OT planning and scheduling is the incorporation of stochasticity to ensure that schedules are either robust or reactive under uncertainty. Uncertainty in the surgical department includes (but is not limited to) patient arrival times, surgical durations, non-elective arrivals, cancellations, staff absence, and OR breakdowns. Deterministic scheduling approaches can become infeasible when faced with even small amounts of schedule deviation and, as such, stochastic approaches are necessary. Whilst robust scheduling strategies can mitigate risks associated with realistic surgical environments, the incorporation of reactive strategies can improve OT performance.

Stuart and Kozan (2012) consider the reactive scheduling of the OT, formulating the SCS Problem (SCSP) as a single-machine scheduling problem. The model includes sequence dependent processing times, due dates, and a lognormal approximation of surgical durations. A bin-packing approach maximises the weighted number of on-time surgeries, considering both elective and non-elective patients. Wang et al. (2015) also apply a machine scheduling approach and formulate the surgical scheduling problem as a no-wait permutation flow-shop. The authors use a predictive-reactive approach for a single surgical suite. Whilst we do not explicitly consider a machine scheduling formulation, the problem in this paper is equivalent to a resource-constrained parallel-machine scheduling problem with identical machines, machine eligibility restrictions, and machine release dates.

Duma and Aringhieri (2015) are the first to present a real-time management model for ORs. The hybrid simulation and optimisation model for the real-time management of ORs incorporates stochastic deviations in surgical duration. Reactions to disruptions include using overtime, cancelling patients, or changing the sequence of remaining patients. A real-time management system increased the number of patients that received their surgery before their surgical due date. The authors have since incorporated non-elective patients (Duma and Aringhieri, 2018) and have evaluated the impact of non-elective capacity reservation policies (Duma and Aringhieri, 2019). Implemented reaction strategies include the resequencing of surgeries and the assignment of overtime. The authors use break-in-moments to ensure that surgeries are not interrupted. Whilst the authors consider multiple non-identical ORs, they do not consider surgeon availability or suitability.

Bruni et al. (2015) consider dynamic scheduling in a multi-OR environment to assign and sequence patients over a one-week planning horizon. The authors incorporate stochasticity in surgical durations and non-elective arrivals, and implement recourse strategies including overtime, swapping, and complete rescheduling. Heydari and Soudi (2016) extend the work of Bruni et al. (2015) to consider both the assignment and sequencing of surgeries under uncertainty, whilst considering the impact of the schedule on the OT and Post Anaesthesia Care Unit (**PACU**). Although Heydari and Soudi (2016) do consider the predictive-reactive rescheduling of multiple ORs, they assume that the ORs are identical, all patients are of the same specialty, and do not consider the assignment of surgical staff.

### 2.1. Innovation

The approach presented in this paper is innovative in a number of ways. Firstly, few authors consider the real-time reactive rescheduling of multiple ORs (Duma and Aringhieri, 2015, 2018, 2019; Heydari and Soudi, 2016). Of the authors that do consider the multi-OR reactive rescheduling problem, only Duma and Aringhieri (2019) consider non-identical ORs. Whilst Duma and Aringhieri (2019) consider non-identical ORs, and allow non-elective patients to be assigned to any of the suitable ORs, they do not allow elective surgeries to be redistributed amongst available ORs during the allocation scheduling phase.

Here, we assume that ORs are non-identical, such that they are suited to certain specialties, but not others. In doing so, we are able to reallocate both elective and non-elective surgeries amongst suitable ORs to reduce the need for overtime. Additionally, we consider the allocation of surgeons to each surgery under surgeon availability and suitability constraints. As such, we do not treat surgeons or ORs as pooled resources, rather individual entities.

Throughout the literature, several authors use break-in-moments to ensure that non-elective patients are treated in a time-sensitive manner, without pre-empting existing surgeries (Wang et al. , 2015; Duma and Aringhieri, 2019). Here, we do not consider break-in-moments, and instead use a hybrid OR dedication policy where some OR capacity is reserved for non-elective surgeries in the original schedule. This policy is flexible in that we are able to reallocate capacity as necessary to reduce the likelihood of overtime. Whilst we do not explicitly consider break-in-moments, we consider do not allow pre-emption of surgeries after anaesthesia has been administered.

In this paper, we provide an interactive OT planning and scheduling methodology that is suited to real-life implementation. This interactive methodology requires analysis of historical data to determine suitable model parameters and assumptions. Here we assume that the master surgical schedule (**MSS**) and surgical case assignments (**SCA**) provide an input for the SCSP, including information on specialty, surgeon, and patient assignments. Through the real-time reactive rescheduling framework presented here, we implement the SCSs and reschedule as necessary to maintain feasibility. We recommend a feedback loop between the SCS and SCA to improve the quality of schedules produced.

Throughout the literature, there is no consideration to scheduling additional elective patients where surgical capacity exceeds pre-scheduled demand or patient cancellations occur. When appropriate surgeons are available and time-to-surgery is sufficient to allow for transport to the hospital and surgical preparation, we schedule additional elective patients. This is in-line with the short-notice scheduling of certain patients at the case study hospital. We also consider the scheduling of non-elective patients. Non-elective patients are those requiring either emergency or urgent surgery. We assume that the scheduler becomes aware of non-elective patients when they arrive at the hospital. The hospital must treat non-elective patients as soon as possible.

The model presented in Section 3 is a bespoke formulation of the reactive SCSP, equivalent to a resource-constrained parallel-machine scheduling problem. In addition to standard resource-constrained parallel-machine scheduling formulations, this problem also requires the consideration of identical machines, machine eligibility restrictions, and machine release dates. The inclusion of identical machines is in reference to the assumption that the duration of a surgery is not dependent on its assigned OR. We use machine eligibility restrictions to ensure that staff only perform surgeries in suitable ORs (i.e. with the correct equipment or features). Machine release dates are used to ensure that patients are not scheduled for surgery in an OR before the OR becomes available. These are innovative additions to the current reactive scheduling literature.

## 3. Model Formulation

We formulate the real-time reactive SCSP is using a MIP approach. We base this on a case study of the surgical department at a large Australian public hospital (cf. Section 5). We provide parameter and variable definitions in the below, before presenting the objective and constraints.

### 3.1. Scalar Parameters

$\overline{H}$: the number of surgeons scheduled for the day.
$\overline{P}$ : the number of patients.
$\overline{S}$ : the number of surgical specialties.
$\overline{R}$ : the number of ORs.
$M$: a sufficiently large number to be used in the Big M constraints.
$\tau$: the schedule start time, in hours since 8am.
$\lambda$: the number of hours that each OR is open for during a standard working day.
$\lambda^*$: the number of hours in a day.

## 3.2. Index Sets

$H$ : the set of surgeons that practice at the hospital. $H = \{1, ..., \bar{H}\}$

$P_S$: the set of elective patients that are scheduled for the day.

$P_U$: the set of patients on the waiting list, that are not included in the current schedule.

$P_N$: the set of non-elective patients to be scheduled for the day.

$P$ : the set of all patients. $P = \{1, ..., \bar{P}\}, P = P_S \cup P_U \cup P_N$

$S$ : the set of surgical specialties). $S = \{1, ..., \bar{S}\}$

$R$ : the set of ORs. $R = \{1, ..., \bar{R}\}$

## 3.3. Indices

$h$: index for surgeon in set $H$.

$p$: index for patient in set (and subsets of) $P$.

$q$: alternative index for patient in set (and subsets of) $P$.

$s$: index for specialty in set $S$.

$r$: index for OR in set $R$.

## 3.4. Vector Parameters

$\alpha_p$: the notice (in hours) an unscheduled patient must receive before surgery, $\forall p \in P_U$.

$\gamma_p$: the arrival time of non-elective patient $p$ (in hours since 8am), $\forall p \in P_N$.

$\mu_p$: the expected duration of patient $p$'s surgery (from start of anaesthesia), $\forall p \in P$.

$\rho_h$: the expected release time of surgeon $h$. This is the time that their current surgery ends, their shift start time, or the time they can arrive at the OR in case of emergency, $\forall h \in H$.

$\kappa_r$: the expected release time of OR $r$, $\forall r \in R$.

$V_p^+$: the setup time (in hours) before patient $p$'s surgery, $\forall p \in P$.

$V_p^-$: the clean-up time (in hours) after patient $p$'s surgery, $\forall p \in P$.

$B_r$: 1 if OR $r$ is working on the day, 0 otherwise, $\forall r \in R$.

$E_{ph}$: 1 if patient $p$ can be treated by surgeon $h$, 0 otherwise, $\forall p \in P, h \in H$.

$I_{ps}$: 1 if patient $p$ is treated by specialty $s$, 0 otherwise, $\forall p \in P, s \in S$.

$T_{rs}$: 1 if OR $r$ is equipped for surgeries by specialty $s$, $\forall r \in R, s \in S$.

## 3.5. Decision Variables

$\epsilon_p$: 1 if patient $p$ is included in the schedule, 0 otherwise, $\forall p \in P$.

$U_{pq}$: 1 if patient $p$'s surgery starts before patient $q$'s surgery ends, 0 otherwise, $\forall p, q \in P$.

$X_{pr}$: 1 if patient $p$ is scheduled for surgery in OR $r$, 0 otherwise, $\forall p \in P, r \in R$.

$Y_{ph}$: 1 if patient $p$ is treated by surgeon $h$, 0 otherwise, $\forall p \in P, h \in H$.

$Z_p$: the expected start time (in hours since 8am) of patient $p$'s surgery, $\forall p \in P$.

$Z_p^*$: the expected end time (in hours since 8am) of patient $p$'s surgery, $\forall p \in P$.

### 3.6. Objective Function

The objective (1) is to maximise total OR utilisation in hours. We calculate utilisation by summing over the time each patient spends in the OR during opening hours, if the patient is included in the schedule (indicated by $\epsilon_p$). Opening hours are 8am through 6pm, represented in the model as time zero and $\lambda$ respectively, where $\lambda = 10$. By considering only OR use during opening hours, the model discourages the allocation of overtime.

Given the nature of the surgical department, it is possible for surgeries to start prior to the standard opening hours, during standard opening hours, and after standard opening hours. This also applies to the completion of surgeries. In particular, this is necessary to ensure that non-elective (both emergency and urgent) patients receive their surgery as soon as possible. As such, consider the contribution to the objective function of each case presented in Figure 1.

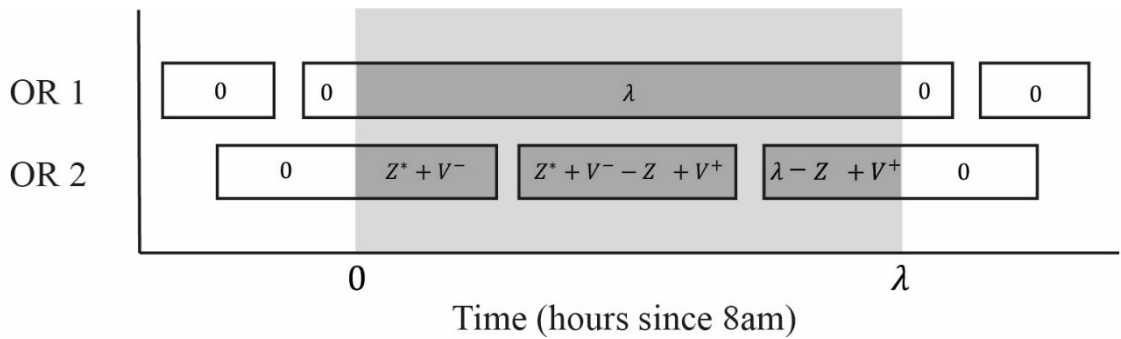

**Figure 1:** Example Gantt chart showing contributions of surgeries to the objective function.

Figure 1 shows the six possible surgery positions with respect to OT opening hours shaded in light grey. Portions of surgeries that contribute to the objective are shaded in dark grey. The portions of surgeries that do not contribute to the objective are unshaded.

To determine the time that a patient $p$ spends in the OR during opening hours, we calculate the difference between the surgery's start and finish times (including setup and clean up), with consideration of the opening hours of the OR. The term $Z_p - V_p^+$ determines the time at which preparation for patient $p$'s surgery begins. The term $Z_p^* + V_p^-$ determines the time at which the clean up after patient $p$'s surgery ends. $Z_p$ and $Z_p^*$ are the expected start and end times of patient $p$'s surgery respectively. $V_p^+$ and $V_p^-$ are the setup and clean up times for patient $p$'s surgery.

The first and last surgeries in OR 1 contribute nothing to the objective as they occur entirely out of OT opening hours. The second surgery in OR 1 contributes a value of $\lambda$ to the objective; however, the portions of time outside opening hours contribute nothing. In OR 2, the first surgery does not contribute to the objective function prior to the OT opening, so its total contribution its end time. The second surgery in OR 2 occurs entirely within opening hours, so its contribution is its duration. The third surgery does not contribute after time $\lambda$.

The contributions presented in Figure 1 can be calculated using the combination of min and max. We use this when forming the objective function (1).

$$\text{Maximise} \quad \sum_{p \in P} \epsilon_p \left( \min(\lambda, \max(Z_p^* + V_p^-, 0)) - \max(0, \min(Z_p - V_p^+, \lambda)) \right) \quad (1)$$

It is possible to linearise this objective function through inclusion of additional decision variables as follows.

$\delta_{jp}^+, \delta_{jp}^-$: non-negative dummy variables used to linearise min and max in the original objective, $\forall j \in \{1, \dots, 4\}, p \in P$

$\varepsilon_{jp}$: 1 if $\delta_{jp}^+$ is greater than zero, zero otherwise, $\forall j \in \{1, \dots, 4\}, p \in P$

$\Omega_p$: the contribution that patient $p$'s surgery makes to the objective function.

$$\text{Maximise} \quad \sum_{p \in P} \Omega_p \quad (2)$$

We use constraints (3) to (5) to ensure that a patient can only contribute to the objective if we schedule them for surgery.

$$\Omega_p \leq \lambda \epsilon_p, \forall p \in P \quad (3)$$

$$\Omega_p \leq \lambda - \delta_{3p}^- - \delta_{4p}^+, \forall p \in P \quad (4)$$

$$\Omega_p \geq \lambda \epsilon_p - \delta_{3p}^- - \delta_{4p}^+, \forall p \in P \quad (5)$$

We use constraint (6) to when linearizing $\max(Z_p^* + V_p^-, 0)$ to $\delta_{1p}^+$.

$$Z_p^* + V_p^- = \delta_{1p}^+ - \delta_{1p}^-, \forall p \in P \quad (6)$$

We use constraint (7) when linearizing $\min(Z_p - V_p^+, \lambda)$ to $\lambda - \delta_{2p}^-$.

$$Z_p - V_p^+ - \lambda = \delta_{2p}^+ - \delta_{2p}^-, \forall p \in P \quad (7)$$

We use constraint (8) when linearizing $\min(\lambda, \max(Z_p^* + V_p^-, 0))$ to $\lambda - \delta_{3p}^-$.

$$\delta_{1p}^+ - \lambda = \delta_{3p}^+ - \delta_{3p}^-, \forall p \in P \quad (8)$$

We use constraint (9) when linearizing $\max(0, \min(Z_p - V_p^+, \lambda))$ to $\delta_{4p}^+$.

$$\lambda - \delta_{2p}^- = \delta_{4p}^+ - \delta_{4p}^-, \forall p \in P \quad (9)$$

Constraints (10) and (11) ensure that only one of $\delta_{jp}^+$ or $\delta_{jp}^-$ can be positive for any combination of $j$ and $p$.get

$$\delta_{jp}^+ \leq M \varepsilon_{jp}, \forall j \in \{1, \dots, 4\}, p \in P \quad (10)$$

$$\delta_{jp}^- \leq M(1 - \varepsilon_{jp}), \forall j \in \{1, \dots, 4\}, p \in P \quad (11)$$

Constraint (12) ensures that the dummy variables introduced are non-negative, whilst constraint (13) ensures $\varepsilon_{jp}$ is binary.

$$\delta_{jp}^+, \delta_{jp}^-, \Omega_p \geq 0 \forall j \in \{1, \ldots, 4\}, p \in P \tag{12}$$

$$\varepsilon_{jp} \in \{0,1\}, \forall j \in \{1, \ldots, 4\}, p \in P \tag{13}$$

### 3.6.1. Alternative Objectives

Whilst the main objective we consider is OT utilisation, hospitals often have performance metrics surrounding other objectives. In this subsection, we define several common alternative objectives: overtime, non-elective time-to-surgery, and number of patients treated. We consider these objectives when evaluating algorithm performance in Section 6.3.

Equation (14) calculates the total overtime used in the day's schedule, where standard operating hours are 0 to $\lambda$ (in hours since 8am).

$$\text{Min} \quad \sum_{p \in P} \epsilon_p (\mu_p + V_p^- + V_p^+) - \min(\lambda, \max(Z_p^* + V_p^-, 0)) + \max(0, \min(Z_p - V_p^+, \lambda))) \tag{14}$$

Equation (15) determines the average non-elective time-to-surgery.

$$\text{Min} \quad \frac{1}{|P_N|} \sum_{p \in P_N} Z_p - \gamma_p \tag{15}$$

We calculate the total number of patients treated using equation (16).

$$\text{Max} \quad \sum_{p \in P} \epsilon_p \tag{16}$$

### 3.7. Constraints

In this subsection, we present the constraints required to produce a well-defined SCS.

Constraint (17) determines the patients scheduled throughout the day. The use of decision variable $\epsilon_p$ allows the planning and scheduling staff to add previously unscheduled patients to the schedule, and ensures that all scheduled and non-elective patients are definitely treated. The set of patients, $P$, includes elective (both scheduled and unscheduled) and non-elective patients.

$$\sum_{r \in R} X_{pr} = \epsilon_p, \forall p \in P \tag{17}$$

Constraint (18) ensures that surgeries are only scheduled in an OR if the OR is working.

$$\sum_{p \in P} X_{pr} \leq \bar{P} B_r, \forall r \in R \tag{18}$$

A patient in OR $r$'s schedule cannot start before OR $r$ is released.

$$Z_p \geq X_{pr} \kappa_r, \forall\, p \in P, r \in R \tag{19}$$

Surgeons cannot treat patients until the surgeon is available. This is either when the surgeon's previous surgery ends, when their shift starts, or when they are able to get to the OR in emergencies.

$$Z_p \geq \sum_{h \in H} Y_{ph} \rho_h, \forall p \in P \tag{20}$$

Figure 2 illustrates how to determine if two surgeries overlap.

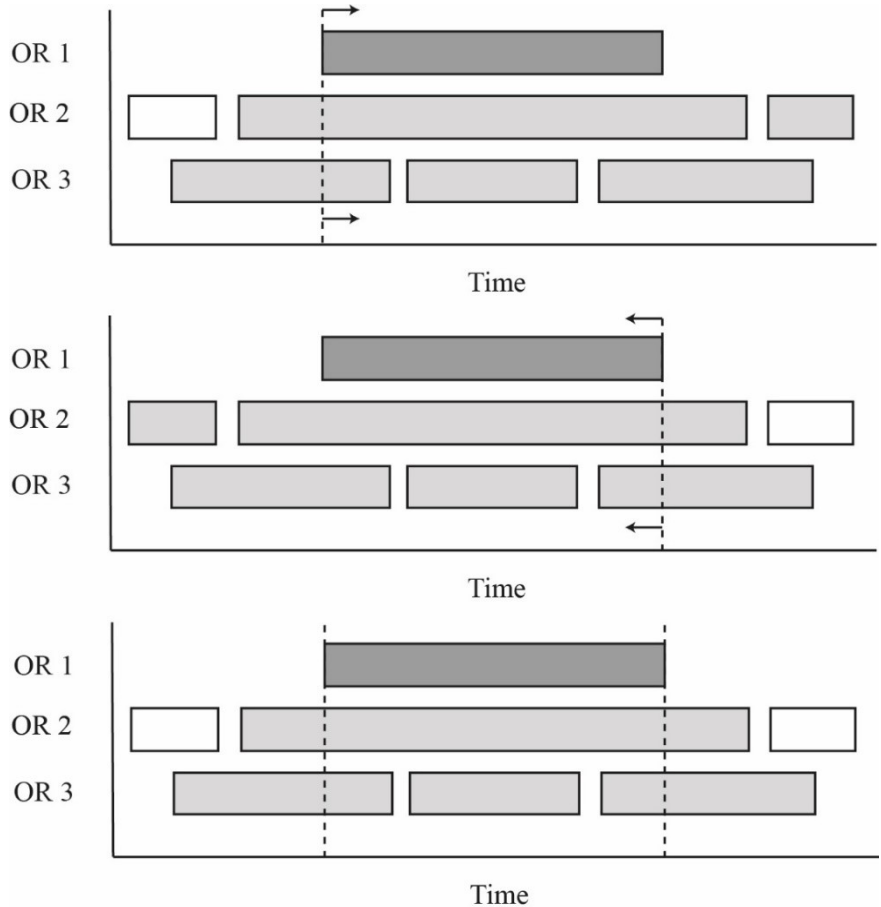

**Figure 2:** Example Gantt chart for determining overlapping surgeries.

In this example, we wish to determine if any surgeries overlap with patient $p$'s surgery shaded in dark grey. First, we select (in light grey) any surgeries that end after patient $p$'s surgery starts. This is equivalent to setting $U_{pq} = 1$. Next, to calculate $U_{qp}$, we select (in light grey) any surgeries that start before patient $p$'s surgery ends. When patient $p$'s surgery begins before patient $q$'s surgery ends, and patient $q$'s surgery starts before patient $p$'s surgery ends, these surgeries overlap.

The binary decision variable $U_{pq}$ is necessary as it represents any patient ordering and determines whether two surgeries overlap. $U_{pq} \times U_{qp} = 1$ if and only if patient $p$ and $q$ have surgeries that overlap. We have linearised this below. Constraints (21) and (22) determine whether patient $p$'s surgery starts before patient $q$'s surgery ends.

$$Z_p - Z_q^* + \lambda^*(2 - \epsilon_p - \epsilon_q) > -\lambda^* U_{pq}, \forall p, q \in P, p \neq q \quad (21)$$

$$Z_q^* - Z_p \geq \lambda^*(U_{pq} - 1), \forall p, q \in P, p \neq q \quad (22)$$

Planning and scheduling staff must not assign a surgeon for overlapping surgeries.

$$Y_{ph} + Y_{qh} \leq 3 - U_{pq} - U_{qp}, \forall p, q \in P, p \neq q, h \in H \quad (23)$$

Constraint (24) ensures that overlapping patients do not receive treatment in the same OR.

$$X_{pr} + X_{qr} \leq 3 - U_{pq} - U_{qp}, \forall p, q \in P, p \neq q, r \in R \quad (24)$$

Constraint (25) calculates the expected finish time of each surgery.

$$Z_p^* = Z_p + \mu_p \epsilon_p, \forall p \in P \quad (25)$$

Constraints (26) and (27) determine the earliest start time of patient $q$, given that patient $p$'s surgery starts before patient $q$'s surgery ends. If the same surgeon (26) or the same OR (27) is used to treat patients $p$ and $q$, then it is necessary to incorporate clean up and setup times between the surgeries.

$$Z_q \geq M(U_{pq} + Y_{ph} + Y_{qh} - 3) + Z_p^* + V_q^+ + V_p^-, \quad (26)$$
$$\forall p, q \in P, p \neq q, h \in H$$

$$Z_q \geq M(U_{pq} + X_{pr} + X_{qr} - 3) + Z_p^* + V_q^+ + V_p^-, \quad (27)$$
$$\forall p, q \in P, p \neq q, r \in R$$

Only certain ORs are suited to certain surgical specialties.

$$X_{pr} I_{ps} \leq T_{rs}, \forall p \in P, r \in R, s \in S \quad (28)$$

Surgeons are qualified or unqualified to perform certain surgeries. Only qualified surgeons can treat patients.

$$Y_{ph} \leq E_{ph}, \forall p \in P, h \in H \quad (29)$$

If a patient is scheduled, then we assign the patient to exactly one surgeon.

$$\sum_{h \in H} Y_{ph} = \epsilon_p, \forall p \in P \quad (30)$$

Elective patients may not start surgery until at least the schedule start time. This may be prior to standard OT opening hours, to account for urgent non-elective treatments.

$$Z_p \geq \tau \epsilon_p, \forall p \in P_S \cup P_U \quad (31)$$

Constraint (32) ensures that previously scheduled elective patients and all non-elective patients are definitely included in the schedule. Previously unscheduled elective patients may be scheduled, but it is not necessary for them to receive treatment on this day.

$$\epsilon_p = 1, \forall p \in P_S \cup P_N \quad (32)$$

Previously unscheduled (add-elective) patients require a certain amount of notice if they are to be treated. This practice is in line with the case study hospital's current requirements for organ transplants.

$$Z_p \geq \epsilon_p(\tau + \alpha_p), \forall p \in P_U \quad (33)$$

Non-elective patients cannot be treated before they arrive at the hospital.

$$Z_p \geq \epsilon_p \gamma_p, \forall p \in P_N \tag{34}$$

Whilst non-elective surgery times are unrestricted, we prohibit the scheduling of add-elective patients if that surgery results in expected overtime.

$$Z_p^* \leq \lambda, \forall p \in P_U \tag{35}$$

Constraint (36) ensures that decision variables are constrained to binary values as appropriate.

$$\varepsilon_p, X_{pr}, Y_{ph}, U_{pq} \in \{0,1\}, \quad \forall p, q \in P, p \neq q, r \in R, h \in H \tag{36}$$

## 4. Solution Methodology

Given the highly uncertain OT environment, schedule disruptions occur frequently. We require a reactive rescheduling framework. Each time a disruption occurs we must update the data and modify the schedule to maintain feasibility. Whilst the model presented in Section 3 could be solved to provide the optimal schedule, the SCSP is NP-hard (Cardoen et al. , 2009). As such, the problem is too time consuming to solve to optimality given the size of the case study instances. This necessitates the use of heuristic techniques that generate good feasible solutions in short amounts of time (less than one second).

In this section, we present a framework for reactive rescheduling of the SCSP (see Figure 3). The main decisions required when reactively rescheduling the OT are decisions regarding schedule frequency and reaction procedure.

Classifications of reactive scheduling frequency include periodic, adaptive, or continuous. Periodic scheduling is when rescheduling occurs after fixed or variable times. Scheduling that occurs after a predetermined amount of schedule deviation is adaptive scheduling. Continuous rescheduling is the strategy whereby rescheduling occurs after a set number of random events.

When rescheduling, many different actions are possible. The extreme responses are those in which either rescheduling of all remaining jobs occurs or there is no response and the system must repair itself. A compromise between these two extremes occurs when the scheduler repairs the system using small schedule changes to maintain feasibility.

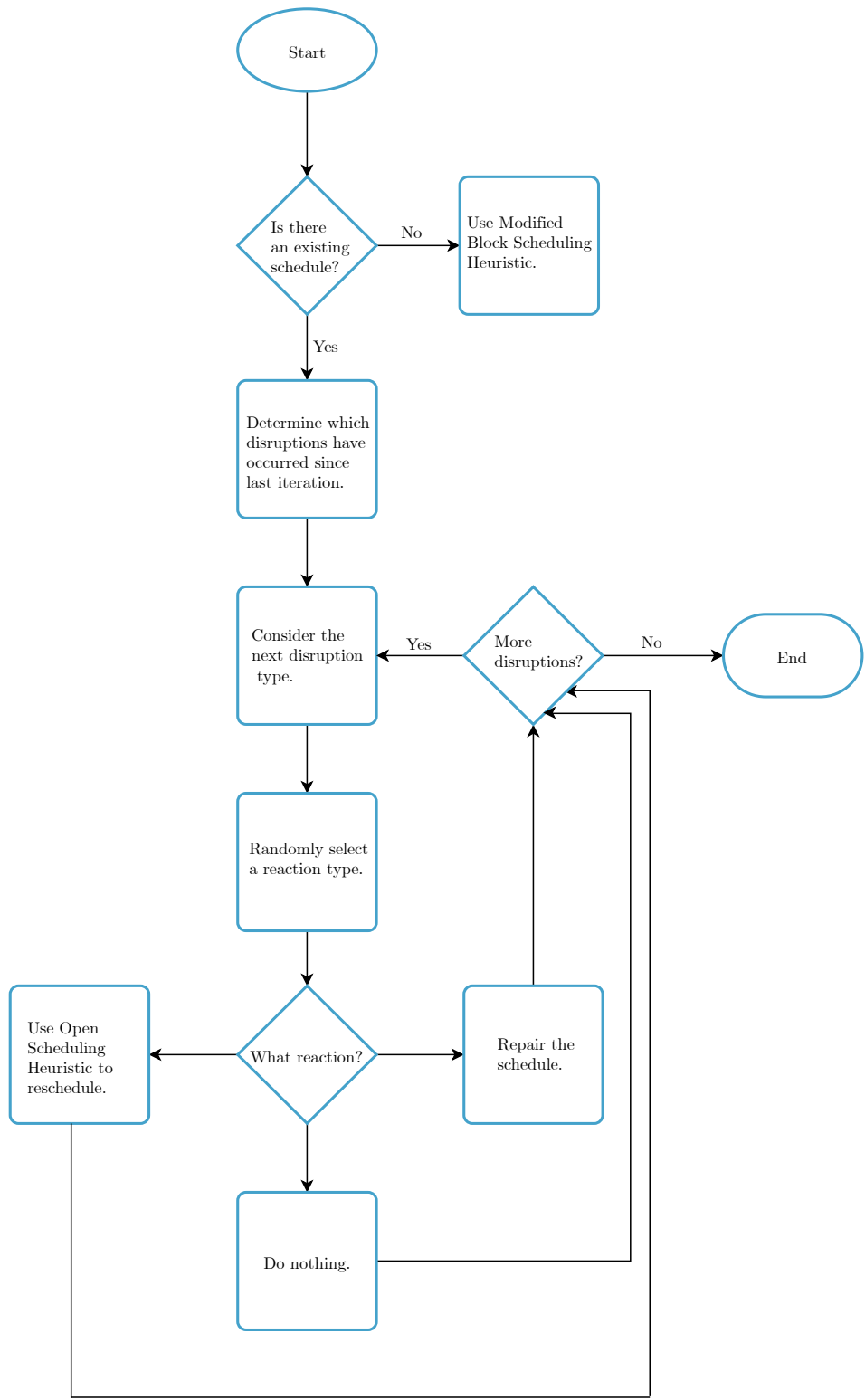

**Figure 3:** Real-time reactive framework.

Here, the entire set of previously scheduled elective and all known non-elective patients are scheduled. There is potential for the day's schedule to include previously unscheduled elective patients where sufficient capacity exists and the patient receives enough notice. Since we assume that surgeries start at the administration of anaesthesia, surgical durations reflect this, and the reactive scheduling techniques prohibit pre-emption of patients once they have been anaesthetised. We also assume that there are sufficient anaesthetists to perform the procedures. This is in-line with the case study hospital.

We investigate all forms of schedule frequency, which in the case of periodic and adaptive frequencies may often lead to a do-nothing reaction to schedule deviations. Within periodic and adaptive frequencies, appropriate levels of responsiveness are determined.

### *4.1. Creating Feasible Schedules*

For solutions to the reactive SCSP to be implementable in real-time, it is necessary to generate solutions within seconds. Here, we present two constructive heuristics for use on the reactive SCSP. A modified block scheduling policy is the basis for the first constructive heuristic. A modified block scheduling policy is a common strategy for OT time allocation where specialties receive time blocks that they can use for their surgical lists. Specialties can then swap blocks as needed. An alternative to a block scheduling is open scheduling. Under an open scheduling policy there are no blocks of time and surgeons can treat patients of any specialty at any time. An open scheduling policy inspired the second constructive heuristic. In each case, we assume that there is a predetermined specialty-OR allocation and a list of required elective surgeries.

#### *4.1.1. Modified Block Scheduling Constructive Heuristic*

This constructive heuristic respects any predefined surgeon assignments as much as possible, before using empty OR blocks to schedule non-elective surgeries. The SCSs generated by the heuristic are in line with a modified block scheduling strategy. Algorithm 1 contains the pseudocode for the modified block scheduling heuristic.

**Algorithm 1:** Modified Block Scheduling Heuristic.
  **FOR** $r$ in the set of working ORs
     Determine the set patients, $P_r$, to be treated in $r$.
     Sort $P_r$ in ascending order of due date.
     **FOR** $p$ in $P_r$
        Append $p$ to OR $r$'s schedule in the earliest possible position.
        Update the surgeon and OR release times.
     **END**
  **END**
  **FOR** $r$ in the set of broken down ORs
     Determine the set patients, $P_r$, to be treated in $r$.
     **FOR** $p$ in $P_r$
        Append $p$ to the schedule in any suitable OR.
        Update the surgeon and OR release times.
     **END**
  **END**
  **FOR** $s$ in the set of non-elective specialties
     **WHILE** there are unscheduled non-elective patients
        **FOR** $r$ in the set of reserved ORs
           Append $p$ to OR $r$'s schedule in the earliest possible position.
        **END**
     **END**
  **END**
  **FOR** $p$ in the set of unscheduled non-elective patients
     Append $p$ to any suitable OR, in the earliest possible position.
  **END**

Algorithm 1 begins by iterating over the set of working ORs, in ascending order. The elective surgery list and existing surgeon allocations indicate the patients treated in the OR. These patients are sorted such that they are treated in order of recommended surgical due date (found by considering both the date of initial surgical request and the patient's urgency category). We append each patient to the OR's schedule in the earliest possible position, by considering any in progress surgeries allocated to the required surgeon, or in progress in the OR.

In some cases, ORs will not be available for use due to equipment failure or other unforeseen events. For each of the unavailable ORs on a given day, we redefine patients scheduled in that OR as urgent non-elective patients. In doing so, these surgeries are not required to follow any pre-existing surgical schedule which is now infeasible due to OR unavailability.

Once we have scheduled or redefined elective patients, the heuristic must assign both surgeons and ORs to all non-elective patients. The algorithm iterates over the set of specialties that have non-elective patients waiting for surgery. We assign each waiting patient to an OR in the set of ORs reserved for that specialty. In each case, we select the non-elective patient who has been waiting longest and append the patient to the schedule with a suitable surgeon in the earliest possible position. The algorithm removes scheduled patients from the set of waiting patients. The process repeats for each specialty until all reserved ORs are used, or there are no non-elective patients of that specialty waiting for surgery. We assign any remaining patients to a suitable OR with an appropriate surgeon. The OR and surgeon assigned to the surgery are selected such that preoperative patient waiting time is minimised.

### 4.1.2. Open Scheduling Constructive Heuristic

The constructive heuristic presented here uses the pre-existing surgeon and specialty assignments as input for surgeon and patient schedules when solving the SCSP. This constructive heuristic uses an open-scheduling strategy and does not necessarily uphold the predefined surgeon-block assignments. Algorithm 2 displays the pseudocode for the Open Scheduling Constructive Heuristic.

For each patient, all feasible surgeon-OR assignments are evaluated. The patient is assigned to the surgeon-OR combination with the earliest possible surgery start time. Whilst this is a greedy constructive heuristic, it tends to perform well on the SCSP under an open scheduling policy. Whilst we could also implement a metaheuristic approach, here it is integral that we produce feasible solutions as quickly as possible. The heuristic maintains feasibility with respect to the MIP SCSP formulation (Section 3). The use of an open scheduling policy incorporates the flexibility required to adjust quickly to disruptions in the OT, whilst maintaining feasibility.

**Algorithm 2:** Open Scheduling Heuristic.
  **FOR** *p* in the set of patients to be scheduled
    Determine patient *p*'s specialty, *s*.
    Determine the set of surgeons who can treat patient *p*.
    Determine the set of ORs that can be used by specialty *s* to treat patient *p*.
    Determine the release time for every possible surgeon-OR combination.
    Select a surgeon-OR combination with the earliest release time.
    Add patient *p* to the schedule, starting at the earliest release time.
    Assign surgeon *h* and OR *r* to patient *p*'s surgery.
  **END**

### 4.2. Real-Time Reaction Strategies

In this subsection, we discuss the reactive strategies used for scheduling under uncertainty. Different strategies are suited to different disruption types. We use D to denote disruption types and define seven main types of disruptions in Algorithm 3. Furthermore, we split strategies again by reaction type: do-nothing (R0), repair (R1), and reschedule (R2).

When generating an initial schedule, we use a modified block scheduling strategy (as per Algorithm 1). When we perform a do-nothing reaction, we simply update the data without making any schedule changes. Note that it is not possible to use a do-nothing approach in response to disruptions of type D2 or D4. When a repair strategy (R1, R1a, or R1b) is used that requires minor rescheduling, or when the reschedule strategy is used, an open scheduling occurs through the constructive heuristic presented in Algorithm 2. We define all intermediate reactions (repair strategies) in Algorithm 3.

When referring to any of the reactive strategies (Algorithm 3), both the disruption type and reaction type are stated. For example, appending an additional patient due to expected OR under-time would be reaction D6R1. Disruption types D6 and D7 occur only in response to disruption types D1 to D5.

When performing the reactions presented in Algorithm 3, it is essential that we respect all constraints (Section 3). For example, if appending a patient to the end of an OR schedule (D1R1, D6R1) it is necessary to ensure sufficient notice is provided to a patient that was previously unscheduled. Thus, at each schedule update, schedules remain feasible (with respect to all constraints).

**Algorithm 3:** Disruption-Reaction Strategy
  **D1.** Non-elective patient arrives
    **R1.** Append patient to the end of an OR schedule.
  **D2.** OR breakdown
    **R1.** Reschedule only the patients assigned to the broken down OR (this must be done by applying an open-scheduling strategy).
  **D3.** Surgery runs under-time
    **R1a.** Shift all remaining surgeries in the OR to an earlier start time (if feasible given surgeon availability).
    **R1b.** Reschedule only this OR and affected surgeon.
  **D4.** Surgery runs over-time
    **R1a.** Shift all remaining surgeries in the OR to a later start time (if feasible given surgeon availability).
    **R1b.** Reschedule only this OR and affected surgeon.
  **D5.** Patient cancellation
    **R1.** Shift remaining surgeries in the OR to an earlier start time (if feasible given surgeon availability).
  **D6.** OR under-time expected
    **R1.** Append an additional patient to the end of this OR's schedule.
  **D7.** OR over-time expected
    **R1.** Reschedule the affected OR and surgeon(s).

The reactions associated with each disruption type are equivalent to local search heuristics. To improve the performance of the reactive rescheduling system, and escape local optima, we randomly select a reaction type based on the disruption that has occurred. Given that some reactions are better suited to particular update frequencies (cf. Section 4.3), we perform parameter tuning to determine the best proportions of each reaction type to use for each combination of disruption and update strategy (cf. Section 6.2).

*4.3. Update Strategies*

There are three main strategies for the frequency of updates when reactively rescheduling: periodic, continuous, or adaptive. Whenever an update strategy triggers an update, we determine the set of disruptions that have occurred since the last reaction. We then perform a reaction according to Algorithm 3.

*4.3.1. Periodic Update*

Periodic updates are updates that occur after set periods. These periods may be constant or may vary throughout the schedule realisation. A periodic update often necessitates a do-nothing reaction. Since do-nothing reactions are not possible for disruptions D2 and D4, periodic updates combined with continuous updates are reactions to disruptions D2 and D4. Algorithm 4 defines the periodic updates.

**Algorithm 4:** Periodic Update Frequency.
------
**UP1.** Update every 0.25 hours. Update continuously with respect to D2 and D4.
**UP2.** Update every 0.5 hours. Update continuously with respect to D2 and D4.
**UP3.** During OT opening hours, update every 0.25 hours. Update continuously with respect to D2 and D4, regardless of time of day.
**UP4.** During OT opening hours, update every 0.5 hours. Update continuously with respect to D2 and D4, regardless of time of day.

*4.3.2. Continuous Update*

Continuous updates ensure that a reaction takes place whenever a schedule disruption occurs. That is, whenever a schedule realisation differs from the original schedule, an update must occur. We define a continuous update strategy, UC, such that we update the schedule whenever a disruption occurs.

*4.3.3. Adaptive Update*

Depending on the situation, periodic updates may occur too infrequently whilst continuous updates may be too frequent. Adaptive updates are a good balance between the two aforementioned update strategies. Rather than performing an update after a set time interval (periodic) or every single disruption (continuous), schedule updates occur after multiple disruptions occur. For example, perhaps three non-elective patients must arrive before a schedule update occurs. Algorithm 5 defines the adaptive updates where updates UA1 through UA5 align with disruption types D1 through D5.

**Algorithm 5:** Adaptive Update Frequency.

**UA1.** Update whenever three or more non-elective patients are waiting to be scheduled.
**UA2.** Update whenever an OR breaks down.
**UA3.** Update whenever a surgery runs under-time by more than 0.5 hours.
**UA4.** Update whenever a surgery runs over-time.
**UA5.** Update whenever a patient is cancelled.

Based on the above adaptive updates, UA2, UA4, and UA5 are equivalent to continuous updates in response to disruptions D2, D4, and D5 respectively. We define an adaptive update strategy (UA) such that an update occurs whenever UA1 through UA5 is triggered.

## 5. Case Study

A case study of a large Australian public hospital is the basis of the work presented in this paper. The hospital performs over 20,000 surgeries, of which around 30% are non-elective patients. Around 360 elective surgery requests occur each week. Administrative staff assign elective surgery requests to one of three urgency categories based on the timeframe within which they should receive their surgery: 30, 90, and 360 days for categories one, two, and three respectively.

Over 110 non-elective surgery requests arise each week. An appropriate surgeon must treat these patients as soon as possible. Whilst non-elective surgeries occur on any day, elective surgeries must occur on weekdays.

At present, the hospital uses a four-week rotating MSS under a modified block scheduling policy, whereby administrative staff (flexibly) assign surgeons or surgical specialties to blocks of OR time. The hospital currently has a ten-hour working day and does not allow the sharing of OR blocks between specialties. The hospital reserves OR blocks throughout the week for non-elective surgery. Surgical staff select elective patients from the waiting list weekly and update this selection frequently to reflect cancellations.

Hospital administrative staff members use general intuition to reschedule surgeries in the event of schedule disruptions. At present, there does not exist a rescheduling policy for maintaining feasibility in real-time. The work in this paper is an improvement on the current reactive scheduling techniques used by administrative staff and improve the utilisation of the OT, without increasing staff overtime.

### *5.1. Data and Parameters*

The data used for computational experiments is a combination of real case study data, and randomly generated data (due to confidentiality) that aligns with historical case study data. The hospital has 21 ORs, 27 specialties, over 100 surgeons, and a waiting list of approximately 2800 elective patients. On average, 113 non-elective surgical patients receive treatment each week.

At present, the surgical department reserves two ORs for non-elective surgeries. As approximately 30% of patients in the OT department are non-elective, the reservation of approximately 10% of ORs is insufficient. The hospital's current strategy has led to the use of overtime and an overflow of non-elective surgeries into ORs assigned to elective surgery. This disruption reduces the efficiency of the ORs, increases surgical overtime, and increases patient cancellations. Along with disruptions caused by the treatment of non-elective patients, stochastic surgical durations, patient cancellations, and machine breakdowns are frequent in the ORs.

To test the model and solution techniques, data is required regarding the waiting list, elective surgery requests, non-elective arrivals, surgical durations, and various cancellations. This historical case study data is available in various forms; however not only is this data confidential, but more information is required when reactively rescheduling the surgical department. Statistical validation of the following data assumptions is presented by Spratt et al. (2018).

In this section, the data required to simulate the hospital environment is summarised. To test the approach presented in Section 4, a waiting list, elective requests, non-elective requests, and cancellations are required. Six problem instances were generated using historical data to verify and validate computational performance [dataset] (Spratt and Kozan, 2018). In the real-life implementation of the approach presented in this paper, the algorithm would use actual hospital data; we cannot present this for confidentiality reasons.

In generating the waiting lists, the number of patients on the waiting list is determined using a Poisson random variable, in-line with historical data. A thinned Poisson process generates the number of patients for every surgeon-specialty-category combination, based on historical averages. Each patient is assigned a surgical duration found by sampling from a lognormal distribution with parameters that vary by surgical specialty. The historical surgical durations are well suited to a lognormal distribution, and informed the choice of parameters. We generate the number of days a patient has been waiting using linear regression based models. These models vary by pre-assigned patient urgency category and consider the impact of a preferred maximum days waiting (as defined in government initiatives to reduce elective surgery waiting time). The parameters to these models were determined through analysis of the historical waiting lists.

We generate both elective and non-elective surgery requests according to a Poisson process with arrival rate dependent on specialty. We sample surgical duration from a lognormal distribution with parameters dependent on surgical specialty. Note that the lognormal duration parameters differ for elective and non-elective patients, and both were determined through analysis of historical data.

In addition to initial waiting lists and surgical requests, it is necessary to include a set of realistic disruptions. Patient based, staff based, and OR based cancellations are considered. It is determined that patients cancel independently of each other and each patient has the same probability of cancelling. Only patient cancellations that occur on the day of surgery are considered. Each OR has a probability of breaking down on each day and OR breakdowns are independent. When an OR breaks down, the OR is unavailable for the entire day.

## 6. Results and Discussion

To verify the real-time reactive SCSP model (Section 3) and the proposed solution methodology (Section 4), computational experiments are performed using datasets based on historical data. We compare the reactive strategies presented in Section 4.2 by running the solution techniques in response to various data realisations and model parameters. In Section 6.2, we perform parameter tuning on a realistic calibration instance to determine the best disruption-reaction probabilities. Following thorough parameter tuning, we use additional test instances for computational experiments in Section 6.3.

Table 1: Summary of test instances.

| Instance | Length of Waiting List | Elective Requests | Non-elective Requests | OR Breakdowns |
|---|---|---|---|---|
| 0 | 2802 | 358 | 120 | 1 |
| 1 | 2759 | 389 | 84 | 0 |
| 2 | 2802 | 352 | 101 | 0 |
| 3 | 2780 | 395 | 100 | 1 |
| 4 | 2744 | 369 | 106 | 4 |
| 5 | 2780 | 373 | 112 | 0 |

These test instances are based on historical data obtained from the case study (cf. Section 5.1). Table 1 contains a summary of the calibration instance (instance zero) and the five test instances used in computational experiments. Each instance is representative of a single week at the case study hospital.

### 6.1. Exact Approach

To provide a computational benchmark, we compare the real-time reactive strategies to a solution obtained using commercial solver Gurobi (Gurobi Optimization, 2018) on a desktop computer with an Intel® Core™ i7 processor at 3.40 GHz with 16.0 GB of RAM. Here, we run the model for each day in the week with the smallest possible subset of data. We do this by only considering previously scheduled patients and non-elective patients, and the set of surgeons, specialties and operating rooms required to treat them. Any non-elective patients that arrive after standard OT hours must wait until the following day. This is as their surgeries cannot contribute to OT utilisation outside of standard OT hours.

When producing these daily schedules, we assume perfect knowledge of future events including cancellations, surgical durations, and non-elective arrivals. Due to the size of the waiting list and computational complexity of the problem, it is not possible to solve the model with short-notice add-on (previously unscheduled) elective patients. We enforced an upper time limit of 1200 seconds for each day. On days where the solution did not converge we report the upper bound on the utilisation.

### 6.2. Parameter Tuning

For each of the update frequencies, we perform parameter tuning for determining the best disruption-reaction probabilities using an iterative approach.

Initially, reaction probabilities are set to a do-nothing approach, the first disruption type (D1) is considered, and the best utilisation is set to zero. We calculate the average utilisation after the metaheuristic is run $n$ times. If the average utilisation is not the best so far, return to the best performing reaction probabilities and move to the next disruption type. If the reaction probabilities have resulted in the best average utilisation, update the best solution values. Randomly perturb the reaction probabilities for the selected disruption type and re-run the real-time reaction strategy.

This approach was repeated until the average quality of solutions (in terms of total OR utilisation) converged. In this instance, 100 iterations were sufficient for UP1, UP2, UA, and UC. After 100 iterations UP3 and UP4 had not converged, however we use the best parameter combination found in that time. Figure 4 displays the results of the parameter tuning in terms of total utilisation (hours) found for each iteration.

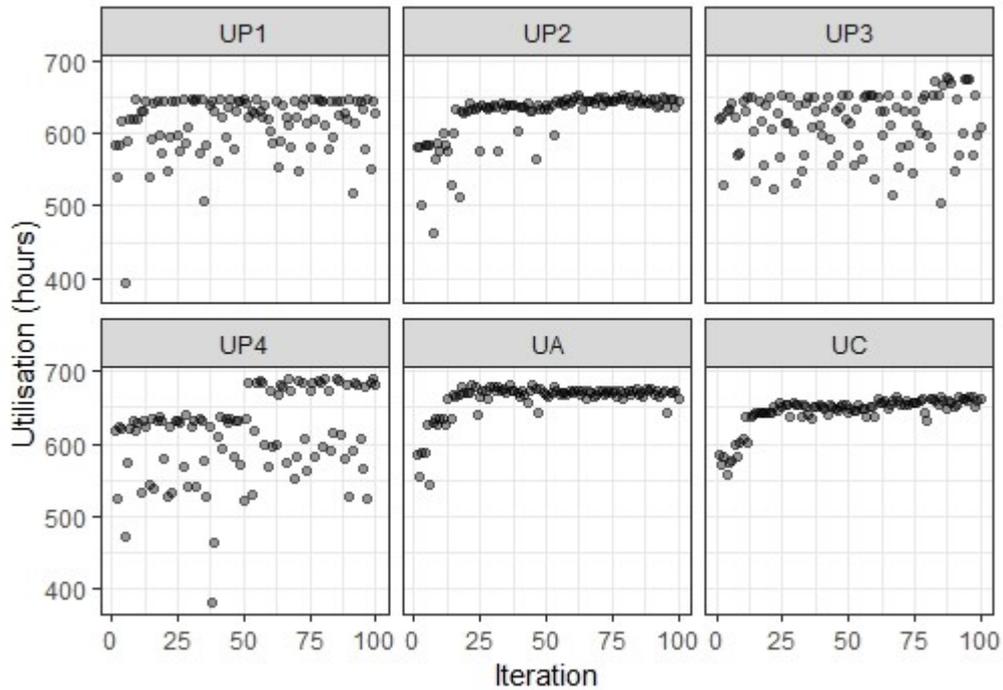

**Figure 4:** Convergence of parameter tuning techniques.

The reaction probabilities provided to the model are used to determine the likelihood of a particular reaction (do nothing, repair, or react) to a particular disruption (D1 to D7). Table 2 shows the results of the parameter tuning. In particular, Table 2 contains the probability of using each reaction type for each combination of update frequency and disruption type.

Based on the parameter tuning, it appears that repair strategies (R1, R1a, and R1b) are often preferred over do-nothing (R0) and reschedule (R2) strategies. Based on the tuned disruption-reaction probabilities (Table 2), we present the best test instance results in Table 3.

**Table 2:** Tuned reaction probabilities for each frequency-disruption combination.

| Disruption | Reaction | Update Frequency | | | | | |
|---|---|---|---|---|---|---|---|
| | | UA | UC | UP1 | UP2 | UP3 | UP4 |
| D1 | R0 | 1.00 | - | 1.00 | 1.00 | 1.00 | 1.00 |
| | R1 | - | 1.00 | - | - | - | - |
| | R2 | - | - | - | - | - | - |
| D2 | R1 | - | - | 1.00 | 1.00 | - | 0.50 |
| | R2 | 1.00 | 1.00 | - | - | 1.00 | 0.50 |
| D3 | R0 | - | 0.50 | - | - | 0.50 | - |
| | R1a | - | - | 1.00 | - | 0.50 | 1.00 |
| | R1b | - | - | - | - | - | - |
| | R2 | 1.00 | 0.50 | - | 1.00 | - | - |
| D4 | R1a | 0.33 | 0.50 | 1.00 | 0.33 | 1.00 | 1.00 |
| | R1b | 0.33 | 0.25 | - | 0.33 | - | - |
| | R2 | 0.33 | 0.25 | - | 0.33 | - | - |
| D5 | R0 | 0.50 | 1.00 | - | - | - | - |
| | R1 | 0.50 | - | 1.00 | 1.00 | 1.00 | 1.00 |
| | R2 | - | - | - | - | - | - |
| D6 | R0 | - | 0.50 | 0.25 | - | - | 0.50 |
| | R1 | 1.00 | 0.50 | 0.50 | 0.50 | 0.50 | 0.50 |
| | R2 | - | - | 0.25 | 0.50 | 0.50 | - |
| D7 | R0 | - | - | 1.00 | - | - | 1.00 |
| | R1a | - | - | - | 1.00 | - | - |
| | R1b | 1.00 | 1.00 | - | - | - | - |
| | R2 | - | - | - | - | 1.00 | - |

**Table 3:** Reactive SCSP Calibration Instance Results

| Update | Utilisation | Runtime (s) | Average Updates | Time per Update (s) |
|---|---|---|---|---|
| UP1 | 645.61 | 12.84 | 852.10 | 0.02 |
| UP2 | 647.11 | 27.26 | 537.90 | 0.05 |
| UP3 | 670.30 | 11.17 | 463.60 | 0.02 |
| **UP4** | **686.17** | **4.27** | **351.20** | **0.01** |
| UA | 672.75 | 16.70 | 504.50 | 0.03 |
| UC | 661.61 | 49.21 | 743.90 | 0.07 |
| Exact | 795.80 | 6022.47 | - | - |

Table 3**Error! Reference source not found.** includes the total OR utilisation over the week, total runtime across a week, the average number of schedule updates that occurred throughout the week, and the average runtime per update. **Error! Reference source not found.** Table 3

also shows the average number of schedules required throughout the week. This number varies based on update strategy and data realisations.

The periodic updates of type four (UP4) performed best out of the update strategies. UP4 required only 4.27s for the entire week's scheduling and was the most efficient of the update strategies. By not allowing updates outside of standard operating theatre times, non-elective patients arriving after the OT closed wait until the next day, thus contributing to the utilisation of the OT. Excluding UP3, which also performed well, all other update strategies would schedule non-elective patients throughout the night, and their surgeries would not contribute to OT utilisation.

UP4 performed well compared to an exact solution provided by Gurobi. The value reported in Table 3 is the upper bound on the optimal solution. Gurobi was unable to converge to an optimal solution within the limit of 1200s for each day. We note that the Gurobi solution used perfect future knowledge, whilst the real-time reactive strategies used predicted values and adapted to disruptions regularly.

Each individual update in the reactive strategies was computationally inexpensive, and as such, it would be reasonable to implement the reaction strategies in real-time. We use the tuned disruption-reaction probabilities (Table 2) when performing computational experiments in Section 6.3.

### 6.3. Computational Performance on Additional Instances

In this section, we perform computational experiments on a number of test instances (cf. Table 1) generated through analysis of historical data [dataset] (Spratt and Kozan, 2018).

We perform the computational experiments using MATLAB® R2017b on a desktop computer with an Intel® Core™ i7 processor at 3.40 GHz with 16.0 GB of RAM. When performing timing, we run the algorithm ten times for each update type. Table 4 shows the average runtime per schedule update. We run the real-time reaction strategies 100 times on the university's High Performance Computing facility to determine the average utilisation, overtime, non-elective time to surgery, and patients treated (cf. Section 3.6).

**Table 4:** Reactive SCSP computational results

| Instance | Update | Utilisation | Overtime | NE Time to Surgery | Patients Treated | Runtime (s) | Update Time (s) |
|---|---|---|---|---|---|---|---|
| 1 | UP1 | 566.97 | 181.73 | **1.72** | **413.20** | 11.14 | 0.01 |
|   | UP2 | 586.75 | 164.03 | 3.59 | 411.82 | 25.36 | 0.05 |
|   | UP3 | 568.78 | **155.96** | 5.21 | 403.23 | 7.55 | 0.02 |
|   | UP4 | 567.75 | 160.23 | 5.05 | 405.78 | 3.78 | 0.01 |
|   | UA  | **600.51** | **154.68** | 2.76 | **413.05** | 24.00 | 0.05 |
|   | UC  | 584.67 | 165.75 | 2.69 | 410.80 | 59.49 | 0.09 |
|   | Exact | 683.64 | 63.48 | 8.03 | 395.00 | 6027.87 | - |
| 2 | UP1 | 599.54 | 243.19 | **1.33** | **415.41** | 12.07 | 0.01 |
|   | UP2 | 616.91 | 195.28 | 4.69 | 401.67 | 27.60 | 0.05 |
|   | UP3 | 646.83 | 179.33 | 3.84 | 407.82 | 7.51 | 0.02 |
|   | UP4 | **662.19** | **168.48** | 3.97 | 409.00 | 4.29 | 0.01 |
|   | UA  | 635.79 | 190.60 | 3.52 | 409.05 | 27.56 | 0.06 |
|   | UC  | 629.30 | 189.58 | 3.45 | 404.23 | 60.72 | 0.09 |
|   | Exact | 762.53 | 141.16 | 8.73 | 402.00 | 6018.20 | - |
| 3 | UP1 | 595.00 | 214.04 | **0.90** | **417.92** | 10.88 | 0.01 |
|   | UP2 | 618.77 | 168.46 | 4.37 | 413.49 | 28.24 | 0.05 |
|   | UP3 | 592.60 | 186.56 | 4.59 | 411.52 | 7.91 | 0.02 |
|   | UP4 | 605.41 | 173.18 | 4.50 | 411.94 | 5.15 | 0.02 |
|   | UA  | 623.48 | 164.61 | 4.07 | 414.60 | 25.38 | 0.06 |
|   | UC  | **630.33** | **156.26** | 4.22 | 412.63 | 60.39 | 0.09 |
|   | Exact | 715.78 | 86.35 | 6.68 | 405.00 | 6051.77 | - |
| 4 | UP1 | 636.26 | 227.11 | **0.72** | **434.99** | 12.13 | 0.01 |
|   | UP2 | 634.20 | 197.05 | 4.52 | 412.11 | 36.85 | 0.07 |
|   | UP3 | 646.35 | 184.76 | 5.54 | 413.85 | 13.41 | 0.03 |
|   | UP4 | **650.46** | **176.19** | 3.70 | 414.55 | 6.28 | 0.02 |
|   | UA  | 634.11 | 197.62 | 3.22 | 413.06 | 33.20 | 0.07 |
|   | UC  | 631.64 | 199.34 | 3.02 | 411.78 | 62.63 | 0.09 |
|   | Exact | 787.03 | 152.92 | 7.56 | 419.00 | 6021.77 | - |
| 5 | UP1 | 581.22 | 226.11 | **1.03** | **438.35** | 11.70 | 0.01 |
|   | UP2 | 614.40 | 187.35 | 4.39 | 429.63 | 30.31 | 0.06 |
|   | UP3 | 625.53 | 154.16 | 2.51 | 425.41 | 8.78 | 0.02 |
|   | UP4 | **638.37** | **142.70** | 2.60 | 426.41 | 5.10 | 0.01 |
|   | UA  | 630.10 | 166.05 | 3.79 | 430.04 | 28.98 | 0.06 |
|   | UC  | 633.36 | 162.20 | 3.60 | 431.16 | 51.97 | 0.07 |
|   | Exact | 730.92 | 77.88 | 9.22 | 419.00 | 6038.46 | - |

Table 4 shows the average results from 100 runs and includes utilisation (in hours), the overtime (in hours), the time-to-surgery of non-elective patients (the time between arrival and treatment in hours), and the number of patients treated. The table also shows the total runtime and runtime per update in seconds.

In Table 4, the best result of each objective is bold, considering only the real-time reaction strategies. We use two-sample t-tests to determine whether objective values are the same for different update types. Where multiple results in the same instance are bolded, two-sample t-tests indicated insufficient evidence ($p > 0.05$) of a significant difference in means.

In three of the five instances, a periodic update (UP4) strategy (updates every half an hour and no updates overnight) produced the schedules with the highest average utilisation. Adaptive updates (UA) and continuous updates (UC) also performed well. The upper bound on the Gurobi solution shows higher utilisation than the real-time reactive strategies, however the Gurobi solution was unable to converge to optimality within 1200 seconds on a minimal daily dataset. In every instance, the reaction strategy that produced the best utilisation was also the reaction strategy that produced the least overtime (as expected). In every instance UP1 (updates every 15 minutes) produced schedules with the lowest non-elective time-to-surgery and highest number of patients treated (i.e. highest number of elective add-ons).

Compared to the exact solution, the real-time reactive strategies did produce schedules with more overtime, however this was due to the fast treatment of non-elective patients. This is reflected in a significant difference in non-elective time-to-surgery. The reactive strategies were able to schedule an average of 10 to 20 additional short-notice elective patients. This was not possible using the exact solution technique due to the size of the waiting list and the computational complexity of the problem.

In the case-study hospital, the average number of patients treated (both elective and non-elective) is approximately 345. Using the approaches presented here, we schedule an average of between 50 and 100 more patients per week (depending on instance). This is dependent on the original MSS SCA used. We produced the MSS SCA implemented for each of the instances by maximising the expected utilisation of the OT, whilst ensuring that the probability of overtime was at most 30%.

In terms of utilisation, historically, the OT at the case study hospital had around 660 hours of utilisation. This was comprised of 515 hours during weekdays, and 145 on the weekend. Whilst the real-time reactive scheduling strategies have a lower total utilisation, we produce schedules with a higher weekday utilisation, and lower weekend utilisation.

Whilst the algorithms here ensure that non-elective patients are treated as soon as possible, in reality, less urgent non-elective surgeries are held until the weekend. In doing so, the hospital is able to attain higher levels of utilisation during weekdays, and lower levels of overtime. As the hospital can postpone non-urgent non-elective surgeries, the hospital has an average of 80 hours of overtime per week. Since this is a decision support tool, we provide recommendations to ensure patients can be treated as quickly as possible, however medical experts can modify the proposed schedules as they see fit.

## 7. Implementation Recommendations

The methodology presented in this paper is suitable for real-life implementation. Hospitals should use an interactive implementation of the real-time reactive SCSP approach. Initially, a MSS-SCA should be determined, assigning patients, surgeons, and specialties to OR blocks. This will assist in staff rostering and the allocation of human resources. At the beginning of each day, the SCSP provides an initial schedule. Throughout the day, the real-time reactive scheduling framework (cf. Section 4.2) should be utilised to ensure schedules remain feasible. At the end of each week, staff should update the data produce a new MSS-SCA to provide OT plans of the upcoming weeks.

Given the requirements of real-time rescheduling, it would be necessary to integrate the software with current hospital IT systems. It is recommended upstream wards (e.g. the emergency department, intensive care unit (**ICU**), and surgical care unit) provide data in real-time to the scheduling software. While in use, the reactive rescheduling framework should provide updates to downstream wards (e.g. PACU and ICU). The current waiting list information system must be accessible from the scheduling software and vice-versa. It is also worth considering integrated implementation of the proposed software with the hospital's bed management software.

At present, planning staff must estimate model parameters and distributions. Planning staff should consider generating a set of real-life based instances to use for the tuning of disruption-reaction probabilities. In future, it may be worth incorporating a real-time approach for the continuous calibration of model parameters. This may result in better prediction of surgical duration based on factors including specialty, consultant, urgency category, and days waiting. This would also enable real-time updates of parameter estimates, ensuring that the software remains up-to-date with minimal human intervention.

## 8. Conclusion

In this paper, we considered the real-time reactive SCSP, with the objective of maximising OT utilisation during standard opening hours. We presented a MIP formulation for the reactive SCSP (Section 3) based on a case study of a large Australian public hospital (Section 5). Due to the complexity of the reactive SCSP, a number of heuristic strategies were presented (Section 4.2) that, whilst aimed at the specific instance of OT scheduling, are applicable to a wide range of machine scheduling problems under uncertainty. We presented constructive heuristics for both the modified block scheduling policy (Section 4.1.1) and the open scheduling policy (Section 4.1.2) as it is integral that hospital administrative staff have access to feasible schedules as quickly as possible. The strategies presented here provide feasible schedule updates in around 0.05 seconds (cf. Table 4).

Whilst there exists reactive rescheduling techniques for the OT department, the consideration of multiple non-identical ORs, surgeon availability and suitability, and long waiting lists is novel (cf. Section 2). We also considered the inclusion of additional elective patients, in-line with the case study hospital. The solution approach is innovative due to its interactive nature. Through this interactive approach, planning and scheduling staff use real-time information in the production of SCSs.

Computational experiments used realistic instances based on a case study of a large Australian public hospital (Section 5). We generated test instances based on historical data to verify and validate the method. Administrative staff should implement the method in real-time with actual hospital data. Results indicate that the reactive rescheduling strategies proposed in this paper maintain schedule feasibility in response to a wide variety of disruptions. The real-time reactive scheduling strategy produces schedules within an ideal range of OT utilization whilst reducing the average time-to-surgery of non-elective patients and increasing the average OT throughput in terms of number of patients.

A number of recommendations were made regarding implementation of the real-time reactive SCSP (cf. Section 7), enabling OT planning and scheduling staff to utilise both the reactive SCSP model (Section 3) and solution methodology (Section 4). These recommendations included the integration with hospital software, the determination of model parameters and distributions, and consideration of both upstream and downstream wards.

**Conflict of Interest**

The authors declare that they have no conflict of interest.